\documentclass[12pt]{article}

\usepackage{a4wide}
\usepackage{latexsym}
\usepackage{amssymb}

\usepackage{mathrsfs}
\usepackage{amsmath,amssymb,amsfonts} 

\usepackage{xspace,theorem,a4,amssymb}

 \usepackage{graphicx,tabularx} 
 \usepackage{color} 
 \usepackage{natbib}
 \usepackage{graphicx}

\usepackage{srcltx}

\RequirePackage[colorlinks,citecolor=blue,urlcolor=blue]{hyperref}
\RequirePackage{graphicx}

\newcommand{\R}{\mathbb R}

\DeclareRobustCommand{\qed}{%
  \ifmmode 
  \else \leavevmode\unskip\penalty9999 \hbox{}\nobreak\hfill
  \fi
  \quad\hbox{\openbox}}
\newcommand{\openbox}{\leavevmode
  \hbox to.77778em{%
  \hfil\vrule
  \vbox to.675em{\hrule width.6em\vfil\hrule}%
  \vrule\hfil}}
\def\1{{\rm 1\mskip-4,4mu l}}

\newtheorem{corollary}{Corollary}

\newtheorem{theorem}{Theorem}

\renewcommand{\thecorollary}{\thesection.\arabic{corollary}}

\theorembodyfont{\rmfamily}
\theorembodyfont{
\rmfamily}

\begin{document}

\title{
{\bf Constructing Bayes Minimax Estimators via Integral Transformations}
}

\author{ Dominique Fourdrinier
\thanks{
Dominique Fourdrinier is Professor,
Universit\'e de Rouen, LITIS EA 4108,
Avenue de l'Universit\'e, BP 12, 76801 Saint-\'Etienne-du-Rouvray, France,
E-mail: Dominique.Fourdrinier@univ-rouen.fr.
}
\and William E. Strawderman
\thanks{
William E. Strawderman is deceased and was Professor, 
Rutgers University, Department of Statistics. 
We gratefully acknowledge Bill Strawderman’s deep and lasting
contributions to statistics, particularly his foundational work in
Bayes minimax theory. His insights have shaped the theoretical
foundations of decision theory, bridging Bayesian and frequentist
approaches, and continue to influence generations of statisticians. 
}
\and Martin T. Wells
\thanks{
Martin T. Wells is Professor,
Cornell University, Department of Statistics and Data Science,
1190 Comstock Hall, Ithaca, NY 14853, USA.
E-mail: mtw1@cornell.edu. 
}
}


\maketitle

\begin{abstract}
The problem of Bayes minimax estimation for the mean of a multivariate normal distribution under quadratic loss has attracted significant attention recently. These estimators have the advantageous property of being admissible, similar to Bayes procedures, while also providing the conservative risk guarantees typical of frequentist methods. This paper demonstrates that Bayes minimax estimators can be derived using integral transformation techniques, specifically through the \( I \)-transform and the Laplace transform, as long as appropriate spherical priors are selected. Several illustrative examples are included to highlight the effectiveness of the proposed approach.

\end{abstract}

\vfill

\noindent AMS 2010 subject classifications: Primary 62C20, 62C15, 62C10; 
secondary 62A15,44A20. 

\noindent Keywords and phrases: Bayes estimate, integral transformations, 
minimax estimate, multivariate normal mean, superharmonic functions. 

\eject

\section{Introduction}\label{section1}
\setcounter{equation}{0}
\setcounter{theorem}{0}
\setcounter{lemma}{0}
\setcounter{corollary}{0}
\setcounter{remark}{0}

The two principal paradigms for estimating the mean of a multivariate normal distribution are the minimax and Bayesian approaches. The minimax framework has been extensively studied, while the Bayesian perspective, although widely applied in certain contexts, is sometimes criticized for its perceived lack of objectivity. In this paper, we investigate the construction of Bayes minimax estimators, which are particularly appealing as they combine the decision-theoretic optimality of Bayes rules with the robustness and conservatism of minimax procedures. As such, these estimators possess desirable properties from both Bayesian and frequentist viewpoints.

In a normal context of dimension $k$, that is, when
$X \sim {\cal N}_k(\theta,I_k)$, Stein \cite{Stein1981} obtains the minimaximity of the
general estimator $\delta$ of the form $\delta(X)=X+\gamma(X)$ through the unbiased estimator of the risk $k+2 \, {\rm div} \gamma(X) + \|\gamma(X)\|^2$.  Thus, if, for every 
$x \in \R^k$, $2 \, {\rm div} \gamma(x) + \|\gamma(x)\|^2 \le 0$,  
then $\delta$ is minimax. 

When (generalized) priors $\pi$ on $\theta$ are considered, 
by the Brown identity formal Bayes estimators are of the form 
$X+\nabla {\rm Log}m(X)$ (where $\nabla$ denotes the gradient and $m$ the 
marginal density), this risk condition becomes, after some calculation, for every 
$x \in \R^k$, $\Delta \sqrt{m(x)}\le 0$ (where $\Delta$ denotes the Laplacian).  Stein \cite{Stein1981} gives examples for generalized Bayes estimates. 

In many cases, the superharmonicity of $\sqrt{m}$ is quite difficult to verify. 
The difficulty has led researchers to consider the superharmonicity of $m$ instead of $\sqrt{m}$. This is reasonable as 
\begin{equation} \label{eq1}
 \forall x\in \R^k \quad 
\Delta\sqrt{m(x)} = \frac{1}{2 \, \sqrt{m(x)}}
\left(\Delta m(x)-\frac{1}{2} \, \frac{\|\nabla m(x)\|^2}{m(x)}
\right) . 
\end{equation}
Therefore, if $\Delta m(x) \leq 0$, then $\Delta\sqrt{m(x)} \leq 0$. A problem
that occurs is that, as shown in Fourdrinier, Strawderman, and Wells
\cite{FourdrinierStrawdermanWells1998AS}, if one uses a proper prior, the
induced marginal cannot be superharmonic. Therefore, if one insists on a 
superharmonic marginal the underlying prior is necessarily improper. However, 
if the function  $\sqrt{m}$ is superharmonic the prior may be proper. 

The results presented in this article focus on the construction of Bayes minimax procedures. In analyzing such estimators, it is useful to distinguish between \textit{construction} and \textit{verification}. Constructing Bayes minimax rules involves identifying a prior that yields the desired minimax properties. This often involves selecting a prior whose marginal distribution satisfies a condition based on superharmonicity, as exemplified in the approaches of \cite{FourdrinierStrawdermanWells1998AS} and \cite{wells2008generalized}.

In contrast, verifying that a given prior leads to a minimax rule is typically more straightforward. Verification generally involves demonstrating that the resulting marginal distribution satisfies an established sufficient condition. For example, \cite{Strawderman1971AMS} employed the Brown identity to show that the shrinkage function associated with the Strawderman prior meets the minimaxity condition introduced by \cite{baranchik1970family}.

The objective of this paper is not to verify the minimaxity of a given Bayes estimator, but rather to address the inverse problem: the construction of Bayes estimators that are minimax. Section 2 presents the formulation of the problem, outlines the necessary setup, and discusses the superharmonicity of \( \sqrt{m} \). Section 3 introduces a general construction method for spherical priors. Section 4 narrows the focus to the subclass of normal-variance mixtures and provides a concrete example. Section 5 concludes with a discussion of the implications of the results.

The findings of this paper complement those of Fourdrinier \textit{et al.} \cite{FourdrinierStrawdermanWells1998AS}, who proposed a related construction method applicable specifically to variance mixture priors. While their approach is somewhat less general, it yields sharper results within that subclass. A broader overview of Bayes minimax estimation can be found in Chapter 3 of Fourdrinier, Strawderman, and Wells \cite{FourdrinierStrawdermanWells2018}.



\section{The Model}\label{section2}
\setcounter{equation}{0}
\setcounter{theorem}{0}
\setcounter{lemma}{0}
\setcounter{corollary}{0}
\setcounter{remark}{0}

Let $X \sim {\cal N}_k(\theta,I_k)$. The central problem of this paper is that of 
constructing spherical Bayes minimax estimators of $\theta$ under 
quadratic loss function $ L(\theta ,\delta )=\|\delta -\theta \|^2.$  Before 
giving the main results, first recall the elements that lead to Bayes minimax 
estimates. \par Assume that $\theta$ is distributed according to a prior
probability measure with density $\pi$. Then the marginal distribution of $X$ 
has density with respect to the Lebesgue measure in $\R^k$ given by 
\begin{eqnarray} \label{marginal}
	 \forall x\in \R^k \quad 
 m(x) = \frac{1}{(2 \, \pi)^{k/2}}
 \int_{\R^k}  {\rm exp} \! 
 \left( \! 
 - \frac{1}{2} \, \|x - \theta\|^2 
   \! \right)
  \pi(\theta) \, d\theta \, .
\end{eqnarray}

	When the prior is assumed to be spherically symmetric, that is,  
$\pi(\theta) = g(\|\theta\|^2)$ for some function $g$, it is characterized 
by its radial density $\lambda$ related to $g$ by, setting $r = \|\theta\|$, 
\begin{equation}\label{prior2}
\lambda(r) 
= 
    \frac{2 \, \pi^{k/2}}{\Gamma(k/2)} \, r^{k-1} \, g(r^2) \, . 
\end{equation}
In that case, the marginal is spherically symmetric as 
well, that is, of the form $m(x)=\ell(\|x\|)$, for some function $\ell$. 
Setting $u = \|x\|$, the exact form of the marginal labeling function $\ell(u)$  
is deduced in Fourdrinier and Wells \cite{FourdrinierWells1996} and equals 
\begin{equation}\label{marginal2}
\ell(u)= 
\frac{\Gamma(k/2)}{2 \, \pi^{k/2}} \, 
\frac{\exp(-u^2/2)}{u^{(k-2)/2}} 
\int ^\infty _0 \frac{\exp(-r^2/2)}{r^{(k-2)/2}} \, I_{(k-2)/2}(u \, r) \, 
\lambda(r) \, dr \, , 
\end{equation}
where $I_\nu$ denotes the modified Bessel function of order $\nu$. The integral 
in (\ref{marginal2}) is related to the $I$-transform which we discuss below. 
 
The Bayes estimate $\delta_\pi(X)$, defined as the minimizer of the Bayes risk 
(see Berger \cite{Berger1985b}, p 17) is given by 
$ \delta_\pi(X) = X + \nabla {\rm Log} m(X)$.  Thus the Bayes estimator 
$\delta_\pi$ is of the form $\delta_\pi(X) = X+\gamma(X)$. 
Stein \cite{Stein1981} shows, under suitable 
integrability conditions, that an unbiased estimator of the risk of 
$\delta_\pi$ equals 
$k + 2 \, {\rm div}\gamma(X)+\|\gamma(X)\|^2 
= 
k + 4 \, {\Delta\sqrt{m(X)} \big/ \sqrt{m(X)}} $ 
where the divergence and the Laplacian are denoted by ${\rm div}$ 
and $\Delta$, respectively. Thus the risk of $\delta_\pi$ equals 
$ E_\theta[\|\delta_\nu-\theta\|^2]
=
k + 4 \, E_\theta [\Delta\sqrt{m(X)} \big/ \sqrt{m(X)}] $ 
where $E_\theta$ is the expectation with respect to 
${\cal N}_k(\theta,I_k)$.  Therefore a sufficient condition for $\delta_\pi$ to 
dominate the usual minimax estimate $\delta_0(X) = X$ is that the square root 
of the marginal density is superharmonic; in this case $\delta_\pi$ will be 
minimax as well. In Appendix \ref{A.1}, we illustrate the difficulty to show the 
superharmonicity of $\sqrt {m}$ for the Strawderman prior 
\cite{Strawderman1971AMS}. 

There is a connection between properness and superharmonicity. 
As is shown in Fourdrinier, Strawderman and Wells \cite{FourdrinierStrawdermanWells1998AS}, the
fact that the prior is proper implies that the marginal cannot be
superharmonic. Consequently, we shall search for priors among those for which
the square root of the marginal density is superharmonic.  It is important to work
with the superharmonicity of the square root of the marginal rather than just
the marginal itself because we wish to investigate the minimaxicity of Bayes rules based on both proper and improper priors. 

We are now in a position to relate the marginal density of $X$ 
to a $I$-transform (see Appendix \ref{A.3}) and its representation $\ell$ in
(\ref{marginal2}). 
Suppose that the prior $\pi$ is a spherically symmetric distribution with radial
density $\lambda$. By definition of the $I$-transform, it follows that (\ref{marginal2}) equals 
\begin{equation}\label{eq33}
	\ell(u) = h(u) \, {\cal I}_{(k-2)/2}[f](u), 
\end{equation}
with 
\begin{eqnarray}\label{eq34}
	h(u) & =& \frac{\Gamma(k/2)}{2 \, \pi^{k/2}} \, u^{(1-k)/2} \, \exp(-u^2/2)
\end{eqnarray}
and 
\begin{eqnarray}\label{eq35}
	f(r) & =& r^{(1-k)/2} \, \exp(-r^2/2) \, \lambda(r) \, . 
\end{eqnarray}
According to (\ref{eq1}), the superharmonicity condition on $\sqrt {m}$ is equivalent to 
\begin{equation}\label{eq9}
 \Delta m(x)- \frac{1}{2} \, \frac{\|\nabla m(x)\|^2}{m(x)} \leq 0  \qquad \forall x\in \R^k . 
\end{equation}
Since $\nabla m(x)=\ell^\prime(\|x\|) \, x / \|x\|$ and $\Delta m(x) = 
\ell'(\|x\|) \, (k-1) / \|x\|+\ell''(\|x\|)$, Condition (\ref{eq9}) becomes, 
through the change of variable $u = \|x\|$, 
\begin{equation}\label{eq31}
\ell' (u) \, \frac{k-1}{u} + \ell''(u) - 
 \frac{1}{2} \, 
 \frac{(\ell'(u))^2}{\ell(u)} 
 \leq 0 \, . 
\end{equation}

The differential inequality (\ref{eq31}) gives rise to a differential inequality
of the $I$-transform of the function $f$ (related to the radial distribution
$\lambda$) given in (\ref{eq33}). We develop this differential inequality in the
next section and use it to characterize classes of minimax priors. 


\newpage

\section{General Spherical Priors}
\setcounter{equation}{0}
\setcounter{theorem}{0}
\setcounter{lemma}{0}
\setcounter{corollary}{0}
\setcounter{remark}{0}

\subsection{The Main Result}\label{subMR}
In this section, we focus on the construction of Bayes minimax rules based on
priors that are spherically symmetric about the origin.  The class of
spherically symmetric priors contains a wide variety of distributions, some of 
which have quite heavy tails. For more on modeling with a normal sampling
distribution and spherical prior, see Angers and Berger \cite{AngersBerger1991CJS}, 
DasGupta and Rubin \cite{DasGuptaRubin1988-procSDTRTIV}, 
Berger and Robert \cite{BergerRobert1990AnnStat}, 
Pericchi and Smith \cite{PericchiSmith1992-JRSS}, 
Fourdrinier and Wells \cite{FourdrinierWells1996} and
Fourdrinier, Strawderman and Wells \cite{FourdrinierStrawdermanWells2018}. 

Theorem \ref{theorem3.1}, 
has two parts. Part 1 expresses Inequality (\ref{eq31}) as previously announced. Part 2 is the main result of this section and gives a construction of a spherically symmetric minimax prior.

\begin{theorem}\label{theorem3.1}
Assume that the prior is spherically
symmetric with the density of the radius $\lambda$ as in \eqref{prior2} and the Bayes risk of the estimate is finite. \par
\begin{enumerate}
\item A sufficient condition for the corresponding spherical Bayes estimator of 
the mean of the multivariate normal distribution to be minimax is
\begin{equation}
\frac{F^{\prime\prime}(u)}{F(u)}
 - 
 \frac{1}{2} \left ( \frac{F^{\prime}(u)}{F(u)} \right )^2
+ 
\frac{F^{\prime}(u)}{F(u)} \left[\frac{k-1}{2} \, \frac{1}{u} - u \right]
+ 
 \frac{(k-1)(7-3k)}{8} \frac{1}{u^2} 
 + 
 \frac{1}{2} u^2 
 - 
 \frac{k+1}{2}
< 0 \label{eq27}
\end{equation}
where $F(u)$ is the $I$-transform of order $(k-2)/2$ of $f$ in \eqref{eq35}. 
\item Let 
\[
 \varphi(u) = \frac{1}{u^2} \sum^{\infty}_{j=0} b_{j} \, u^{j} 
\]
be any non positive generalized series. The function F defined by
\begin{equation} \label{eq28} 
F(u) = (c_{1}\,  z_{1}(u) + c_{2} \, z_{2}(u))^{2} \, u^{(k-1)/2}\,  \exp (u^{2}/2) \, , 
\end{equation} 
where $z_{1}$ and $z_{2}$ are two linearly independent solutions of
the second order differential equation
\begin{equation} \label{eq29} 
z^{\prime \prime} + \frac{k-1}{u} \, z^{\prime} - \frac{1}{2} \, \varphi(u) \, z = 0
\end{equation}
and $c_{1}$ and $c_{2}$ are arbitrary constant, satisfies (\ref{eq27}).

Provided the so obtained function $F$ is the $I$-transform of order $(k-2)/2$ of a function,
the radial density of the prior is given, for every $r>0$, by
\begin{equation}
\lambda(r) = r^{(k-1)/2} \exp(r^2/2){\cal I}^{-1}_{(k-2)/2}[F](r),
\label{eq30}
\end{equation}
where ${\cal I}^{-1}_{(k-1)/2}$ is the inverse $I$-transform of order $(k-1)/2$. 

If, in addition, the function $\lambda$ is integrable, then the resulting proper
Bayes estimator is minimax. 
\end{enumerate}
\end{theorem}

\noindent{\tt Proof}
According to (\ref{eq33}), we have $\ell(u) = h(u) \, F(u)$ where $h(u)$ is 
given by (\ref{eq34}). Hence (\ref{eq31}) becomes 
\begin{eqnarray}\label{eq36}
[p(u)F(u)+F'(u)] \, \frac{k-1}{u} 
&+&(p^2(u)+p'(u))F(u)+2p(u)F'(u)+F''(u)\nonumber\\
&-& 
 \frac{(p(u)F(u)+F'(u))^2}{2 \, F(u)}
 \leq 0 \, , 
\end{eqnarray}
where $F(u) = {\cal I}_{(k-2)/2}[f](u)$ and $p(u) = (1-k) / 2 \, u^{-1} - u$.
Upon defining $R(u) = F'(u)/F(u)$ Inequality (\ref{eq36}) can be written as
\begin{equation}
[p(u)+R(u)] \frac{k-1}{u}  + \frac{1}{2} \, (p(u) + R(u))^2 + (p'(u) + R'(u)) \leq 0 
\label{eq37}
\end{equation}
which, setting $y(u) = p(u) + R(u)$, becomes
\begin{equation}
y'(u) + \frac{k-1}{u} \, y(u) + \frac{1}{2} \, y^2(u) \leq 0 \, . \label{eq38}
\end{equation}

In order to solve this first order non linear differential inequality define a function
$\varphi \leq 0$ and solve the generalized Ricatti equation
\begin{equation}
y'(u) + \frac{k-1}{u} \, y(u) + \frac{1}{2} \, y^2(u) = \varphi(u) \, . \label{eq39}
\end{equation}
This equation may be solved using the classical method for generalized Ricatti
equations.  First, set $y (u) = \lambda \, v(u) \, u^{1-k}$, 
$r(u) = 1 / 2 \, \lambda \, u^{1-k}$, $s(u) = - 1 / \lambda \, \varphi (u)\,  u^{k-1}$ 
for some $\lambda$. Then solving (\ref{eq39}) reduces to
\begin{equation}
v^{\prime} + r(u) \, v^{2} + s(u) = 0 \, .    \label{eq40}
\end{equation}
Now applying Section 16.515 of Gradshteyn and Ryzhik \cite{GradshteynRyzhik1994book}, 
the general solution of (\ref{eq38}) is 
\[
 v(u) = \frac{1}{r(u} \, 
 \frac{c_{1} \, z^{\prime}_{1} (u) + c_{2} \, z^{\prime}_{2}(u)}{c_{1} \, z_{1} (u) + c_{2} \, z_{2} (u)} \, , 
\]
where $z_{1}$ and $z_{2}$ are independent solutions of the second order differential equation 
\begin{equation}
z^{\prime \prime} - \frac{r^{\prime} (u)}{r(u)} \, z^{\prime} + r(u) \, s(u) \, z = 0
\label{eq41}
\end{equation}
and $c_{1}$ and $c_{2}$ are arbitrary constants.  Upon evaluation of $r$ and 
$s$, (\ref{eq41}) reduces to (\ref{eq29}). 

The hypothesis on $\varphi$ allows us
to construct solutions even when (\ref{eq29}) has a singularity at 0 
(see Krasnov, Kiss\'elev and Markarenko \cite{KrasnovKisselevMarkarenko1981-book}, p. 167). 
Transforming back to (\ref{eq39}), it follows that the general solution of that equation is 
\[
 y(u) = 2 \, \frac{c_{1} \, z^{\prime}_{1}(u) + c_{2} \, z^{\prime}_{2} (u)}
 {c_{1} \, z_{1}(u) + c_{2} \, z_{2} (u)}
  \, .
\]
Hence, by definition of $y(u)$, 
\[
\frac{F^{\prime}(u)}{ F(u)}
 = 
 2 \, \frac{c_{1} z^{\prime}_{1}(u) + c_{2} \, z^{\prime}_{2} (u)}{c_{1} \, z_{1}(u) + c_{2} \, z_{2}(u)}
 + 
 \frac{k-1}{2} \, \frac{1}{u} + u \, .
\]
Upon solving this elementary first-order equation, it follows that
\[
 F(u) = (c_{1} \, z_{1} (u) + c_{2} \, z_{2}(u))^{2} \, u^{(k-1)/2} \, e^{u^{2}/2} \, ,
\]
which proves (\ref{eq28}). 

In order to deduce the form of a (generalized) Bayes minimax prior density, recall that
$F(u)={\cal I}_{(k-2)/2} [f](u)$.
Therefore, ${\cal I}^{-1}_{(k-2)/2}[F](r)=f(r)$ and applying (\ref{eq35})
and solving for $\lambda(r)$ give the result in (\ref{eq29}) and complete
the proof of Part 2. \hfill$\Box$

\vspace{.3cm}

\noindent {\bf Remark}

According to (\ref{eq28}), the function $F(u)$ is nonnegative, so that the radial
density $\lambda$ in (\ref{eq30}) is also nonnegative.


\subsection{An Example of a General Spherical Prior}
Choose $\varphi (u) = - 2 \, b / u^2$ with $b \geq    0$.  Then (\ref{eq29}) is
\begin{equation}
z^{\prime \prime} + \frac{k-1}{u} \, z^{\prime} + \frac{b}{u^{2}} \, z = 0 \, , 
\label{eq42}
\end{equation}
which is a Bessel type equation, that is, 
$z^{\prime \prime} + p(u) \, z^{\prime} + q(u) \, z = 0$ 
(see Appendix \ref{A.2}) where the coefficients $p(u)$ and $q(u)$ can be 
represented as generalized series of the form 
\[
 p(u) = \frac{1}{u} \sum^{\infty}_{j=0} a_{j} \, u^{j} 
 \quad {\rm and} \quad 
 q(u) = \frac{1}{u^{2}} \sum^{\infty}_{j=0} b_{j} \, u^{j}
\]
with $a_{0} = k-1, \ b_{0} = b, \ a_{j} = b_{j} = 0 \ {\rm for} \ j \geq 1.$  
According to Krasnov, Kiss\'elev and Markarenko 
\cite{KrasnovKisselevMarkarenko1981-book},  
the solutions of (\ref{eq42}) are of the form 
$z = u^{\rho} \sum_{j=0}^{\infty} c_{j} \, u^{j}$ (with $c_{0} \neq 0$) 
where $\rho$ is a solution of the determinant equation 
\begin{equation}
\rho \, (\rho - 1) + a_{0} \, \rho + b_{0} = 0 \, .  \label{eq43}
\end{equation}
In this specific case the roots of (\ref{eq43}) are 
\[
 \rho_{1} = \frac{2 - k - \sqrt{(k-2)^{2} - 4 \, b}}{2}
 \quad {\rm and} \quad 
 \rho_{2} =\frac{2 - k + \sqrt{(k-2)^{2} - 4 \, b}}{2} 
\]
for $b \leq (k-2)^{2} / 4$.  When the difference of the roots, 
$\rho_{2} - \rho_{1}$, is not integer or zero, two independent solutions of
(\ref{eq42}) are
\[
 z_{1} (u) = u^{\rho_{1}} \sum^{\infty}_{j=0} c_{j} \, u^{j} 
 \qquad {\rm and} \qquad  
 z_{2} (u) = u^{\rho_{2}} \sum^{\infty}_{j=0} c_{j} \, u^{j} \, , 
\]
where the coefficient $c_{j}$ can be determined by plugging the series 
solutions into (\ref{eq42}), which yields
\begin{equation}
\sum^{\infty}_{j=0} \left[(j + \rho) (j + \rho + k - 2) + b\right] 
c_{j} \, u^{j + \rho -2} =0    \qquad      \forall u \geq 0,  \label{eq44} 
\end{equation}
for $\rho$ equals $\rho_{1}$ or $\rho_{2}$.  The solution to (\ref{eq44}) is 
$c_{j} = 0$ for all $j \geq 1$ and $c_{0} \neq 0$ arbitrary since the first term 
of the series is null by (\ref{eq43}).  Hence the solutions to (\ref{eq42}) are 
\[
 z_{1} (u) = A_{1} u^{\rho_{1}} \qquad {\rm and} \qquad z_{2} (u) = A_{2} u^{\rho_{2}}
\]
for $A_{1}$ and $A_{2}$ arbitrary constants
\footnote{As a side remark, notice that it can be shown that, when $\rho_{2} - 
\rho_{1}$ is an integer or zero, an additional solution of (\ref{eq42}) is 
$z_{3} (u) = A \, z_{1}(u) \log u + B \, z_{2}(u)$.
}.
Hence, the ${\cal I}$-transform of $r \rightarrow r^{(1-k)/2} \exp (r^{2}/2)
\lambda (r)$ is
\[
 F(u) = (A_{1} \, u^{\rho_{1}} + A_{2} \, u^{\rho_{2}})^{2} \, u^{(k-1)/2} \, \exp (u^{2}/2).
\]

As in Theorem \ref{theorem3.1}, the determination of the radial density of the 
prior which yields a minimax Bayes rules reduces to linear combination of terms 
of the form 
\begin{equation}
r^{(k-1)/2} \exp(r^{2}/2) {\cal I}^{-1}_{(k-2)/2} [u^{\gamma} \exp(u^{2}/2)] (r). 
\label{eq45}
\end{equation}
\noindent From Oberhettinger \cite{Oberhettinger1972book} Formula 5.13, 
on the inversion of Hankel transformation, 
\[
 {\cal H}_{\nu} [u^{\eta-1/2} \, e^{- \alpha u^{2}}](r) 
 \propto 
 r^{-1/2} \, \exp \! \left(- \frac{r^2}{8 \, \alpha}\right) 
 M_{\eta /2, \nu/2} \! \left(- \frac{r^2}{4 \, \alpha}\right) , 
\]
with ${\cal R}e(\eta + \nu) > - 1$ and ${\cal R}e(\alpha) > 0$ where
$M_{x, \mu}$ is the Whittaker function connected with the confluent
hypergeometric function (see Appendix \ref{A.2}) by 
\[
 M_{x, \mu} (z) = \exp \left(- \frac{z}{2}\right) \, z^{\mu + 1/2} \, 
 _{1}F_{1} \left(\mu + \frac{1}{2} - x; 1 + 2 \mu ; z\right).
\]
Therefore, using (\ref{eq26}), the expression in (\ref{eq45}) is proportional to
\begin{equation}
 r^{(k-2)/2} \exp \! \left(\frac{r^{2}}{4}\right) 
 M_{\gamma/2 + 1/4, (k-2)/4}\left(\frac{r^{2}}{2}\right).  \label{eq46}
\end{equation}
with $\gamma + (k + 1) / 2 > 0$.  If $b$ is chosen such that 
$(k-2)^{2} / 4 - 1 < b \leq (k-2)^{2} / 4$ all the regularity conditions above 
are valid.  Thus the construction of Theorem \ref{theorem3.1} yields a radial prior 
which gives rise to a minimax Bayes estimator. 

It is interesting to note the striking similarity between (\ref{eq46}) and the
Strawderman prior in (\ref{eq22}).  However, note that the expression
in (\ref{eq46}) does not have a finite integral.  Perhaps choosing an
appropriate linear combination of solutions will yield integrability, 
although we were unable to find one. 

\vspace{.3cm}

\noindent {\bf Remark:} The second-order differential equation in (\ref{eq29}) 
has various other solutions depending on the choice of $\varphi$.  A convenient 
class of solutions can be determined by the generalized Bessel's differential 
equation $ u^{2} z^{\prime \prime} + (1 - 2 \alpha) u z^{\prime} + [(\beta 
\gamma u^{\gamma})^{2} + \alpha^{2} - v^{2} \gamma^{2}] z = 0.$  The solution to 
this equation is $ z(u) = u^{2} Z_{v} (\beta u^{\gamma})$ where $Z_{v}$ is a 
solution of the Bessel's differential equation $ u^{2} z^{\prime \prime} + u 
y^{\prime} + (u^{2} - v^{2}) z = 0.$  Special solutions $Z_{v}$ of this equation 
are the Bessel, Neumann and Hankel functions $J_\nu, Y_\nu, H^{(1)}_\nu, 
H^{(2)}_\nu$, respectively. 


\section{The Variance Mixtures of Normals Case}
\subsection{The Main Result}
\setcounter{equation}{0}
\setcounter{theorem}{0}
\setcounter{lemma}{0}
\setcounter{corollary}{0}
\setcounter{remark}{0}

In this section, we apply the general theory of the previous section and give a
method to construct a family of priors which give rise to the proper Bayes
minimax rules for scale mixtures of normals.

Notice that, in the case where the prior is a scale mixture of normals, we have
that 
\[
 \lambda (r)= \int^\infty_0 h(v) \lambda_v (r) dv
\]
where 
\[
\lambda_v(r) = 
 \frac{2^{1-k/2}}{\Gamma (k/2)} \, \frac{r^{k-1}}{v^{k/2}} \exp \! \left(\! -\frac{r^2}{2 \, v}\right)
\]
is the normal radial density.  Hence the marginal labeling function has the representation
\begin{eqnarray*}
\ell(u) 
 &=& 
 g(u) \, {\cal I}_{(k-2)/2} \! \left[\int^\infty _0 \lambda_v(r) \, h(v)
 \exp\!\left(\!-\frac{r^2}{2}\right) r^{(1-k)/2} \, dv \right]\! (u) \\
 &=& 
 g(u) \int^\infty _0 {\cal I}_{(k-2)/2} [\lambda_v (r) h(v) \, \exp\!\left(\!-\frac{r^2}{2}\right) 
 r^{(1-k)/2}](u) \, dv 
 \\
 &=& 
 \frac{2^{1-k/2}}{\Gamma (k/2)} \, g(u) 
 \int^\infty _0 \! h(v) \int^\infty _0 \! 
  \sqrt{r \, u} \, I_{(k-2)/2}(ru) \, 
  \frac{1}{v^{k/2}} \, 
  r^{(k-1)/2} \, \exp\!\left(\!-\frac{r^2}{2 \, v}\right) 
 \\
 & & \hspace{10cm} \exp\!\left(\!-\frac{r^2}{2}\right) dr \, dv \, . 
\end{eqnarray*}
Then it follows that 
\begin{eqnarray}
 \ell(u) 
 &=& 
 \frac{2^{1-k/2}}{\Gamma (k/2)} \, g(u) \, 
 \sqrt{u} \int^\infty _0 h(v) \frac{1}{v^{k/2}} \,  
 \int^\infty _0 r^{k/2} \exp \!\left(\! - \left(1 + \frac{1}{v}\right) 
 \frac{r^2}{2}\right) 
  \\
 & & \hspace{8.5cm} I_{(k-2)/2} (r \, u) \, dr \, dv  \nonumber\\
 &=& 
 \frac{2^{1-k/2}}{\Gamma (k/2)} \, g(u) \, u^{(k-1)/2} 
 \int^\infty _0 h(v) \frac{1}{v^{k/2}} 
 \left(1 + \frac{1}{v}\right)^{-k/2} 
 \exp \! \left(\frac{u^2 \, v}{2 \, (1+v)}\right) dv  \nonumber\\
 &=& 
 \frac{1}{(2 \, \pi)^{k/2}}
 \int^\infty _0 \frac{h(v)}{(v+1)^{k/2}} 
 \exp \! \left(\frac{- u^2}{2 \, (1+v)}\right) 
 dv \, .  \label{mixmarg}
\end{eqnarray}

\newpage

\begin{corollary}\label{corollary4.1}
 Assume that the prior is a variance mixture of normal distributions with
mixing density $h$.\par
\begin{enumerate}
\item A sufficient condition for the corresponding spherical (generalized) Bayes
estimator of the mean of the normal distribution to be minimax is
\begin{equation}
 \frac{G^\prime(s)}{G(s)} - 2 \, \frac{G^{\prime\prime}(s)}{G^\prime(s)} 
 \le \frac{k}{s} \qquad \forall s>0
\label{eq7}
\end{equation}
where $G$ is the Laplace transform of 
\[
 f : t\rightarrow t^{k/2-2} h \! \left(\frac{1-t}{t}\right){\bf 1}_{]0,1[}(t) 
\]
 ( ${\bf 1}_{]0,1[}$ being the indicator function of the interval
 $]0,1[$ ).\label{item1}
\item For any continuous function $\varphi$ such that, $\varphi (s) \le k / s$ 
for every $s>0$, the function $G$ defined, for any $s>0$, by 
\begin{equation}
G(s)=\left(\int_b^s {\rm exp} \left(-\frac{1}{2} \int_a^t \varphi(u) \, du\right) dt
\right)^2, \label{eq8}
\end{equation}
where $a$ and $b$ are two real constants, satisfies (\ref{eq7}).
Provided that the function $G$ obtained is the Laplace transform of a
positive function defined on $]0,\infty[$ and vanishes outside $]0,1[$, 
the mixing density $h$ is given for every $v>0$, by
\[
 h(v)=(v+1)^{k/2-2} {\cal L}^{-1}[G] \left(\frac{1}{v+1}\right) , 
\]
where ${\cal L}^{-1}$ is the inverse Laplace transform. If, in
addition, the function $h$ is integrable, then the resulting proper Bayes
estimator is minimax. \label{item2}
\end{enumerate}
\end{corollary}

\noindent{\tt Proof}
Making the identification $G(u^{2}/2)=\ell(u)$, where $\ell(u)$ is defined in 
(\ref{marginal2}), it follows from some tedious calculations that (\ref{eq27}) 
reduces to (\ref{eq7}). This proves Part 1.

It is clear that any function $\varphi$ which satisfies 
\[
 \varphi (s) = \frac{G^\prime(s)}{G(s)} - 2 \, \frac{G^{\prime\prime}(s)}{G^\prime(s)} 
 \qquad \forall s>0
\]
yields a solution of (\ref{eq7}) provided that $\varphi(s) \le k / s$. Straightforward calculations lead to
\[
\left({\rm Log} \ 2\left(\sqrt {G(s)}\right)^\prime\right)^\prime = 
 - \frac{\varphi(s)}{2} \qquad \forall s> 0
\]
and two integrations give the following expression for the Laplace transform $G$: 
\[
 G(s) = \alpha\left(\int_b^s {\rm exp} \left(-\frac{1}{2}\int_a^t
 \varphi(u) \, du\right) dt \right)^2 
 \qquad \forall s>0
\]
where $a\in \R$ and $b\in \R$.

Finally, inverting the Laplace transform $G$, we obtain
\[
 t^{k/2-2} h\left(\frac{1-t}{t}\right){\bf 1}_{]0,1[}(t)={\cal L}^{-1}[G](t)
\qquad \forall t\in ]0,1[
\]
that is 
\[
 h(v) = (v+1)^{k/2-2} {\cal L}^{-1}[G]\left(\frac{1}{v+1} \right) 
 \qquad \forall v>0.
\]
This proves Part \ref{item2} of the corollary. 
\hfill$\Box$


\subsection{An Example of a Variance Mixture of Normal Prior}
Corollary \ref{corollary4.1} seems difficult to use directly, due to
the requirement that the function $G$ is a Laplace transform of a positive
function on $]0,1[$.

If this is not the case, an obvious guess would be $\varphi(s) = k / s$ which
gives $G(s) = \alpha\left(s^{1-k/2}+\beta\right)^2$ with $\alpha>0$ and $\beta
\in \R$. When $\beta =0$, we would have
\[
 f(t) = {\cal L}^{-1}[G](t) = \alpha{\cal L}^{-1}[s^{2-k}](t) = 
 \alpha \, \frac{t^{k-3}}{(k-3)!} 
 \qquad \forall t\in ]0,1[.
\]
Hence the corresponding mixing density $h$ equals 
\[
 h(v) = \alpha \, (v+1)^{k/2-2} \, f\!\left(\frac{1}{v+1}\right) 
 = \alpha \, \frac{(v+1)^{1-k/2}}{(k-3)!}    
 \qquad \forall v> 0
\]
which is a special case of the form of Strawderman's prior. 

However, on closer inspection, the analysis is faulty, as $s^{2-k}$ is the
Laplace transform of $t^{k-3} / (k-3)!$ on all $]0,\infty[$. Since the
Laplace transform is unique, $s^{2-k}$ is not the Laplace transform of $t^{k-3} 
/ (k-3)! {\bf 1}_{]0,1[}(t)$. It is quite interesting that the faulty analysis
does give a prior which is proper and for which the resulting Bayes estimation
is minimax. 

It is easy to see from another perspective that $\varphi(s) = k / s$ cannot
correspond to the mixing distribution for a proper prior. If so, the
resulting Bayes estimator, $\delta$, would have risk identically to $k$,
the risk of the usual estimator $\delta_0=X$. In this case, since $\delta$
cannot be equal to $X$ (an unbiased estimator cannot be a proper Bayes prior in the normal
case), $1 / 2 \, (\delta_0 + \delta)$ would have a Bayes risk smaller than
$\delta$. Hence $\delta$ cannot be Bayes, and the ``boundary'' case
$\varphi(s) = k / s$ cannot correspond to a proper prior.

But first we indicate that it is possible to use the theorem and its proof
directly, although we bypass guessing at a suitable $\varphi$. Instead, we first
guess a suitable inverse Laplace transform and show that its Laplace transform
satisfies (\ref{eq7}), and that the function $h$ that results is positive and
integrable. 

\subsubsection{Example 1}

Take ${\cal L}^{-1}[G](t)= t^n {\bf 1}_{]0,1[}(t)$ where $n$ is a 
nonnegative integer. Then an easy calculation gives 
\begin{eqnarray*}
 G(s) & =& 
 \frac{n!}{s^{n+1}} - e^{-s} \sum_{j=0}^n \frac{n!}{j! \, s^{n-j+1}} \\
      & =& 
 \frac{n!}{s^{n+1}} \, e^{-s} \sum_{j=n+1}^\infty \frac{s^j}{j!}  \, , 
\end{eqnarray*}
hence 
\begin{eqnarray*}
 G^\prime(s) = - \frac{(n+1)!}{s^{n+2}} \, e^{-s} \sum_{j=n+2}^\infty \frac{s^j}{j!}  
\end{eqnarray*}
and 
\begin{eqnarray*}
 G^{\prime\prime}(s) = \frac{(n+2)!}{s^{n+3}} \, e^{-s} \sum_{j=n+3}^\infty \frac{s^j}{j!} \, . 
\end{eqnarray*}

Therefore the condition 
\[
 \varphi(s) = \frac{G^\prime(s)}{G(s)} - 2 \, \frac{G^{\prime\prime}(s)}{G^\prime(s)} 
 \le \frac{k}{s} 
\]
becomes 
\[
 - \frac{n+1}{s} \, \frac{\sum_{j=n+2}^\infty \frac{s^j}{j!}}{\sum_{j=n+1}^\infty \frac{s^j}{j!}} 
 + 
 2 \, \frac{n+2}{s} \, \frac{\sum_{j=n+3}^\infty \frac{s^j}{j!}}{\sum_{j=n+2}^\infty \frac{s^j}{j!}} 
 \le \frac{k}{s} 
\]
or, equivalently, 
\begin{equation}
 2 \, (n+2) \left[
 s \, \frac{\sum_{j=n+2}^\infty \frac{s^j}{(j+1)!}}{\sum_{j=n+2}^\infty \frac{s^j}{j!}} 
 \right]
 - (n+1) \left[
 s \, \frac{\sum_{j=n+1}^\infty \frac{s^j}{(j+1)!}}{\sum_{j=n+1}^\infty \frac{s^j}{j!}} 
 \right]
 \le k \, . 
\label{eq11}
\end{equation}

Now we note that both terms in brackets in (\ref{eq11}) are less than 1. This
follows since each term in the numerator is contained in the denominator. 
Additionally, note that the first bracketed term is smaller than the second
bracketed term. In fact, each of the two terms is of the form $s$ times
$E[ 1 / 1+X]$ where $X$ is a Poisson random variable truncated at $n+2$ and $n+1$, 
respectively. Since the Poisson truncated at $n+2$ is stochastically larger than 
that truncated at $n+1$ the claim follows immediately.  Then (\ref{eq11}) is
true as soon as $2 \, (n+2)-(n+1)\leq k$ or equivalently $n\leq k-3$.\par It remains
to find conditions under which the corresponding density is integrable. Since 
the mixing density is given by
\begin{eqnarray*}
 h(v) 
 & =& 
 (v+1)^{{k / 2}-2} \, {\cal L}^{-1}[G]\left(\frac{1}{v+1}\right) \\
 & =& 
 (v+1)^{{k / 2}-2-n}     \, .
\end{eqnarray*}
This is integrable when $n > k / 2 - 1$. Hence, the resulting proper Bayes prior
procedure is minimax when $k / 2 - 1 < n \leq k-3$.  When $k=5$, this reduces to
$1.5 < n \leq 2$ or $n=2$ which corresponds to the Strawderman prior to
$a={1 / 2}$. 

While this proof is perhaps less efficient than Strawderman's original result,
it demonstrates the potential utility of the main theorem of this section. 

\subsubsection{Example 2}

The next example is reminiscent of the generalized beta mixtures of
Gaussians, a family of global-local shrinkage priors commonly used in
high-dimensional Bayesian regression (see Armagan, Dunson, and Clyde \cite{armagan2011generalized}). These priors are designed to offer
adaptive shrinkage, robustness, and computational tractability, while
maintaining proper posterior behavior even when the number of
predictors is large. This family provides a flexible and efficient
framework for sparse estimation, unifying and extending many existing
priors. Notably, it encompasses several important priors as special
cases, including the horseshoe, normal-exponential gamma,
Strawderman-Berger, normal-gamma, and normal-beta prime priors. 

In particular, take 
\[
 {\cal L}^{-1}[G](t)= t^{\alpha-1} (1-t)^{\beta - 1} (1-\sigma t)^{-\gamma} {\bf 1}_{]0,1[}(t) \, , 
\]
with $\sigma > 0$. Then 
\[
 G(s) = \int_0^1 f_s (t) \, dt 
\]
with 
\[
 f_s (t) = t^{\alpha-1} (1-t)^{\beta - 1} (1-\sigma t)^{-\gamma} \, e^{-st} 
\]
so that 
\begin{align*}
- \, G'(s) & = \int_0^1 t \, f_s (t)dt\\
& = \int_0^1 t^\alpha \, (1-t)^{\beta-1} \, (1-\sigma t)^{-\gamma} \, e^{-st} \, dt\\
& = \left[t^\alpha \, (1-t)^{\beta - 1} \, (1-\sigma t)^{-\gamma}  \left(-\frac{1}{s}\right) e^{-st}\right]^1_0\\
&\qquad \qquad + \frac1s \int_0^t \frac{d}{dt} (t^\alpha \, (1-t)^{\beta-1} \,  (1-\sigma t)^{-\gamma}) \, e^{-st} \, dt\\
& = \frac1s \int_0^1 \left[\frac{\alpha}{t} \,  t \, f_s(t) - 
\frac{(\beta - 1)}{1 - t} \, 
t \, f_s (t) + \frac{\gamma \, \sigma}{1-\sigma t} \, t \, f_s (t)\right] dt\\
& = \frac{1}{s} \int_0^1 \left[\frac{\alpha}{t} - 
\frac{\beta - 1}{1-t} + \frac{\gamma \, \sigma}{1-\sigma t} \right] t \, f_s (t) \, dt\\
& = \frac1s \int_0^1 \left[\alpha - (\beta - 1) \, \frac{t}{1-t} + 
\gamma \, \sigma \frac{t}{1-\sigma t} \right] f_s(t) \, dt \, . 
\end{align*}
Expressing $G^{\prime\prime}$ the same way gives 
\begin{eqnarray}
 \varphi(s) &=& \frac{G^\prime(s)}{G(s)} - 2 \, \frac{G^{\prime\prime}(s)}{G^\prime(s)} 
 \nonumber \\
 &=& 
 -\frac{1}{s} E_s^* \left[\alpha - (\beta - 1) \frac{T}{1-T} + \gamma \sigma \frac{T}{1-\sigma T}\right] 
 \label{E*} \\
 &+& 
 \frac{2}{s} E_s^{**} \left[(\alpha + 1) - (\beta - 1) \frac{T}{1-T} + \gamma \sigma \frac{T}{1-\sigma T}\right] 
 \label{E**} 
\end{eqnarray}
where $E^*_s$ is the expectation with respect to a density proportional to
$f_s(t)$ and $E_s^{**}$ is the expectation with respect to a density
proportional to $t \, f_s(t)$. 

\medskip
\noindent {\it Case 1: Assume $\gamma < 0$}

\medskip
Note that $t \, f_s(t)$ has a monotonic likelihood ratio with respect to $f_s(t)$ and the integrand in (\ref{E*}) and (\ref{E**}) is decreasing. Then
$$
E_s^{**} \! \left[ (\alpha + 1) - (\beta - 1) \, \frac{T}{1-T} + \gamma \, \sigma \,  \frac{T}{1-\sigma T}\right] 
\leq 
E_s^* \! \left[ (\alpha + 1) - (\beta - 1) \frac{T}{1-T} + \gamma \, \sigma \frac{T}{1-\sigma T}\right] . 
$$
Therefore
\begin{align*}
\varphi (s) & \leq 
\frac{1}{s} \, 
E^*_s\!\left[ \left\{  2 \, (\alpha + 1 ) - \alpha \right\} - (\beta - 1) \frac{T}{1-T} + 
\gamma \, \sigma \frac{T}{1-\sigma T}\right] \\
& = 
\frac1s \left\{\alpha + 2 + 
E_s^*\!\left[ (1-\beta) \, \frac{T}{1-T} + \gamma \, T \, \frac{T}{1-\sigma T}\right]\right\} \\
& \leq \frac{\alpha  + 2}{s} \, , 
\end{align*}
since the integrand term is negative.

\medskip
\noindent {\it Case 2: Assume $\gamma > 0$}

Similarly
\begin{align*}
\varphi (s) & \leq \frac1s \left\{
\alpha + 2  + E_s^*\!\left[ (1-\beta) \, \frac{T}{1-T}\right] + 2 \, E_s^{**}\! \left[\gamma \, \sigma \, \frac{T}{1-\sigma T}\right]\right. \\
&\qquad \left. - E_s^*\!\left[\gamma \, \sigma \, \frac{T}{1-\sigma \, T}\right]
\right\} . 
\end{align*}
Hence
\begin{align*}
\varphi (s) & \leq \frac1s \left\{\alpha +2 + E_s^*\!\left[(1-\beta) \,  \frac{T}{1-T}\right] + 2 \, E_s^{**}\!\left[\gamma \, \sigma \frac{T}
{1-\sigma \, T}\right]\right\}\\
& \leq \frac1s \left\{ \alpha +2 + \frac{2 \, \gamma \, \sigma}{1-\sigma}\right\} , 
\end{align*}
since $\gamma \, \sigma \, T / (1-\sigma \, T)$  is increasing and 
$(1-\beta) \, T / (1-\sigma \, T) \leq 0$.

Therefore the minimaxity condition $\varphi (s) \leq k / s $ 
is expressed as follows. 

Case 1: $$\frac{\alpha +2}{s} \leq \frac{k}{s} \Leftrightarrow \alpha + 2 \leq k$$

Case 2: $$\alpha + 2  + 2 \frac{\gamma \, \sigma}{1 - \sigma}\leq k \, . $$


\section{Concluding Remarks}

For the class of spherical priors, we have shown that the superharmonicity of the square root of the marginal density leads to a differential inequality that, in many cases, can be solved using integral transformation techniques to yield a Bayes minimax rule. The methods developed in this paper provide a useful framework for constructing default priors for Bayesian analysts who seek estimators with favorable frequentist properties. We hope that the proposed technique can be extended to construct Bayes minimax estimators in a broad range of settings. While it would be desirable to apply our approach to construct proper Bayes minimax rules, the associated computations are currently intractable. The results of this paper complement those of Fourdrinier, Strawderman, and Wells \cite{FourdrinierStrawdermanWells1998AS}, who proposed an alternative construction method that is specifically applicable to priors defined as variance mixtures of normal distributions.


\setcounter{section}{1}
\renewcommand{\thesection}{\Alph{section}}
\section*{Appendix}\label{Appendix}

\setcounter{equation}{0}
\setcounter{subsection}{0}
\setcounter{corollary}{0}
\renewcommand{\thecorollary}{\Alph{corollary}}
\setcounter{example}{0}
\renewcommand{\theexample}{\Alph{example}}
\setcounter{lemma}{0}
\setcounter{remark}{0}
\setcounter{theorem}{0}

\subsection{On the condition $\Delta\sqrt{m}\le 0$ 
for the Strawderman prior} \label{A.1}

Note that it turns out that the condition $\Delta\sqrt{m}\le 0$ in \eqref{eq1} 
can be difficult to verify. For instance, Strawderman's prior \cite{Strawderman1971AMS} was
obtained through a two-stage hierarchical model, that is conditional on 
$\lambda$ $(0 <\lambda \le 1)$, $\theta$ is normal with mean 0 and covariance 
matrix $\lambda^{-1} \, (1-\lambda) \, I_k$ and the unconditional density of 
$\lambda$ with respect to the Lebesgue measure is given by $(1-a) \, 
\lambda^{-a}$ for any $a$ such that $0 \le a < 1$.  Hence, this gives rise to 
the prior density 
\begin{eqnarray}\label{eq19}
	\pi(\theta) 
	&=& 
	(2 \, \pi)^{-k/2} \int^1 _0 \left(\frac{\lambda}{1 - \lambda}\right)^{k/2}
	\exp \! \left(\frac{- \lambda \, \|\theta\|^2}{2 \, (1 - \lambda)}\right)
	(1-a) \, \lambda^{-a} \, d \lambda  \nonumber\\
	&=& 
	(1-a) \, (2 \, \pi)^{-k/2} \int^\infty _0 t^{k/2-a} \, (1+t)^{a-2} \, 
	\exp \! \left(- \frac{t \, \|\theta\|^2}{2}\right) dt \, .  
\end{eqnarray}
Therefore $\pi(\cdot)$ is a spherically symmetric density and hence by 
(\ref{eq19}) has radial density
\begin{equation}\label{eq21}
	\lambda(r) = \frac{2 \, (1-a)}{2^{k/2} \Gamma({k / 2})} \, r^{k-1} 
	\int^\infty _0 t^{k/2-a} (1+t)^{a-2} \, 
	\exp \! \left(- \frac{t \, r^2}{2}\right) dt \, . 
\end{equation}

In Magnus, Oberhettinger and Soni \cite{MagnusOberhettingerSoni1966}, it is
shown that the Whittaker function, $W_{x,\mu} (z)$, is related to the integral
in (\ref{eq21}). In particular, they provide (p.313) the following integral 
representation of the Whittaker function, that is, 
\begin{eqnarray*}
	\Gamma  \! \left(\frac{1}{2} + \mu - x\right) W_{x, \mu} (z) = 
	\exp \! \left(\frac{- z}{2}\right) z^{\mu + 1 / 2} 
	\int^\infty _0 \exp(-z \, t) \, t^{\mu -  x - 1 / 2} \, (1+t)^{\mu + x - 1 / 2} 
	\, dt.
\end{eqnarray*}
Therefore, applying this representation, it follows that (\ref{eq21}) equals
\begin{equation}\label{eq22}
	\lambda(r)= \frac{1-a}{2^{k/4-1}} \left(\frac{k}{2} - a \right) r^{(k-2)/2} 
	\exp \! \left(\frac{r^2}{4}\right) 
	W_{a-1-k/4, (k-2)/4} \! \left(\frac{r^2}{2}\right). 
\end{equation}
Hence, the Strawderman prior has a convenient representation as a
mixture of Whittaker functions. 

According to (\ref{marginal2}), the marginal can be derived through
straightforward but tedious calculations. It is given by 
\begin{eqnarray*}
	\ell(u) = \frac{1-a}{2 \, \pi^{k/2}} \, 
	\frac{\Gamma(k/2-a+1)}{\Gamma(k/2-a+2)}  \,  {}_1F_1(k/2-a+1;k/2-a+2;-u^2/2) \, 
	, 
\end{eqnarray*}
where we used the intergral representation of the confluent hypergeometric 
${}_1F_1 (\alpha;\beta;z)$. Then the minimaxity condition $\Delta\sqrt{m(x)}\le 
0$ reduces to 
\begin{eqnarray*}
	&& -\frac{k/2 + a - 3}{(k/2 - a + 2)^2} \, \frac{\|x\|^2}{2} 
	{}_1F{}_1^2(1;{k/2} - a + 3;\frac{\|x\|^2}{2})
	+ \\
	&& \hspace{5cm} 2 \, \frac{2 - a -\|x\|^2 / 2}{k/2 - a + 2} 
	{}_1F{}_1(1;{k/2} - a + 3;\frac{\|x\|^2}{2}) - 2 
	\le 0 \, , 
\end{eqnarray*}
where $a>3-k/2$. It is easy to check that this inequality is satisfied for $x$ in the neighborhood of the origin and at infinity when $k/2 + a - 3 \geq 0$. 
However, it is considered considerably more difficult to establish the inequality directly
for all $x$. 

\subsection{Some special functions} \label{A.2}

The modified Bessel functions \( I_\nu(x) \) and \( K_\nu(x) \) are 
independent solutions to the modified Bessel differential equation
\[
x^2 \frac{d^2 y}{dx^2} + x \frac{dy}{dx} - (x^2 + \nu^2) y = 0,
\]
where $\nu$ and $z$ can be arbitrary complex valued.  Further, recall their series representation
\[
I_\nu(x) = \sum_{k=0}^{\infty} \frac{1}{k! \, \Gamma(k + \nu + 1)} 
\left( \frac{x}{2} \right)^{2k + \nu} \text{and}\; \; K_\nu(x) = \frac{\pi}{2} \frac{I_{-\nu}(x) - I_\nu(x)}{\sin(\nu \pi)} 
\]
for non-integer \( \nu \). To define \( K_n(x) \) for integer \( n \), one takes 
the limit
\[
K_n(x) = \lim_{\nu \to n} K_\nu(x) = \lim_{\nu \to n} \frac{\pi}{2} \frac{I_{-\nu}(x) - I_\nu(x)}{\sin(\nu \pi)} \, . 
\]

The Bessel function has a number of integral representations under certain conditions.  We use the following
\[
I_\nu(x) = \frac{\left( \frac{x}{2} \right)^\nu}{\Gamma\left( \nu + \frac{1}{2} \right) \sqrt{\pi}} \int_{-1}^{1} \frac{e^{x t}}{(1 + t^2)^{\nu + \frac{1}{2}}} \, dt 
 \, . 
\]
Note that this expression has an interesting probabilistic interpretation. Let \(\mathbf{U} \sim \mathrm{Unif}(S^{n-1})\) and \(\mathbf{v} \in \mathbb{R}^n\) be fixed unit vectors.   Then for \( n = 2 \, \nu + 1 \), the modified Bessel function of the first kind can be written as
\[
I_\nu(x) = \frac{\Gamma(\nu + 1)}{(x/2)^\nu} \, \mathbb{E}\left[ e^{x \, \mathbf{v} \cdot \mathbf{U}} \right] \, . 
\]

	The generalized hypergeometric function is defined as:
\[
{}_pF_q(a_1, \dots, a_p; b_1, \dots, b_q; z) = \sum_{k=0}^{\infty}
\frac{(a_1)_k \cdots (a_p)_k}{(b_1)_k \cdots (b_q)_k}
\frac{z^k}{k!}
\]
where \( (a)_k \) is the Pochhammer symbol (rising factorial) \( (a)_k = a(a+1)(a+2) \cdots (a+k-1)\), for $(a)_0 = 1$. When $p = q = 1$, 
it can be shown that the generalized hypergeometric function $_1F_1(a; b; z) $ is a solution of the Kummer differential equation
\[
z \frac{d^2y}{dz^2} + (b - z) \frac{dy}{dz} - a y = 0 \, . 
\]

	A comprehensive study of the various Bessel and hypergeometric
functions can be found in Magnus et al. \cite{MagnusOberhettingerSoni1966}. 

\subsection{Integral Transformations} \label{A.3}

A definition of the Hankel transform can be found, for example, in Koh and
Zemanian \cite{KohZemanian1969-SIAM} and is as follows. For any $\nu \in \R$ and 
$a>0$ let ${\cal J}_{\nu,a}$ be the function space generated by the kernel 
$\sqrt {xy}\ J_\nu(xy)$, where $J_\nu(.)$ is the Bessel function of the first 
kind, $x>0$, and $y$ is a complex number in the strip $\Omega = \{y:|\hbox {Im}\ 
y|<a,\ y\not= 0\ \hbox {or a negative number}\}$. The Hankel transform ${\cal 
	H}_\nu$ is defined on the dual space ${\cal J}^*_{\nu,a}$ as 
\[
{\cal H}_\nu[f](y) = \int ^\infty _0 f(x) \, \sqrt {xy} \, J_\nu(xy) dx \, , 
\]
where $y\in \Omega$, $f$ is locally integrable and 
\[
\int ^\infty _0 |f(x) \, e^{ax} \, x^{a+1/2}| \, dx < \infty.
\]
Now to every $f\in {\cal J}_{\nu,a}^*$ there exists a number $\sigma_f$ 
(possibly infinite) such that $f\in {\cal J}_{\nu,b}^* $ if $b\le \sigma_f$ and 
$f\notin {\cal J}_{\nu,b}$ if $b > \sigma_f$. Hence $f\in {\cal 
	J}_{\nu,\sigma_f}$. 

Now, for any $f\in {\cal J}_{\nu,\sigma_f}$, define the $I$-transform as
\begin{equation}\label{eq24}
	F(y) = {\cal I}_\nu[f](y) = \int ^\infty _0 f(x) \, \sqrt {xy} \, I_\nu (xy) \, dx
\end{equation}
where $I_\nu(.)$ is the modified Bessel function of the first kind. The 
real-valued $I$-transform is defined where $y$ is restricted to the strip 
$\chi_f = \{s : 0 < s < \sigma_f\}$. The inversion formula for the $I$-transform 
is given by 
\begin{equation}\label{eq25}
	f(x) = {\cal I}_\nu^{-1}[F](x) =
	\lim_{r  \rightarrow \infty} \frac{1}{i \, \pi}
	\int ^{\sigma+ir}_{\sigma-ir} F(y) \, \sqrt{xy} \, K_\nu (xy) \, dy 
\end{equation}
for $\nu\ge -1/2,\ \sigma\in \chi_f$ and where $K_\nu(.)$ is the modified
Bessel function of the second kind of order $\nu$.

The statement of the inverse formula is given by Koh and Zemanian in 
\cite{KohZemanian1969-SIAM} and may be proved using the ideas similar to the 
proof of Theorem 6.6 provided by Zemanian in \cite{Zemanian1968book}. 
It is important to note that the inversion holds only for $\sigma \in \chi_f$, 
therefore $\sigma =0$ is not an admissible value. If one is going to apply 
(\ref{eq25}) it is necessary to use a contour integration technique not just a 
simple change of variables onto the real axis.\par However the inverse 
$I$-transform is also related to the Hankel transform since the Bessel 
function, $J_\nu$, and the modified Bessel function, $I_\nu$, are connected 
by the identity 
\[
I_\nu (z) = \exp(-\pi \nu i/2) \, J_\nu(iz) \quad  {\rm for} \quad -\pi<\arg z \leq \pi/2
\]
(see formula 8.406 of Gradshteyn and Ryzhik \cite{GradshteynRyzhik1994book}). 
Then it follows that 
\[
{\cal I}^{-1}_{\nu}[F](x) = \exp(-i\pi (\nu/2 - 3/4)) \, {\cal H}^{-1}_{\nu}[\hat F](x) 
\]
where $\hat F(y) = F(-iy)$.  As the Hankel transform is self reciprocal, i.e.
$ {\cal H}^{-1}_{\nu}[.]= {\cal H}_{\nu}[.]$,
the previous identity becomes 
\begin{equation}\label{eq26}
	{\cal I}^{-1}_{\nu}[F](x) = \exp(-i\pi (\nu/2 - 3/4)) \ {\cal H}_{\nu}[\hat F](x).
\end{equation}
Both of those relationships are useful when using tables of Hankel's 
transforms. 

The $I$-transform is also related to, the more well known, $K$-transform 
\[
{\cal K}_\nu[g](y) = \int ^\infty _0 g(x) \, \sqrt {xy} \, K_\nu (xy) \, dx \, . 
\]
If $G(y)={\cal K}_\nu[g](y)$ the inversion formula is given by 
\[
{\cal K}_\nu^{-1}[G](x) =
\lim_{r  \rightarrow \infty} \frac{1}{i \, \pi} 
\int ^{\sigma+ir}_{\sigma-ir} G(y) \, \sqrt{xy} \, I_\nu (xy) \, dy,
\]
where $\nu\ge -1/2$ and $\sigma\in \chi_g$. As 
\[
K_\nu(z) = \frac{\pi}{2 \, {\rm sin}(\pi \, \nu)} \, [I{_\nu}(z) - I_\nu(z)]
\]
for non integer $\nu$, it follows that 
\[
{\cal I}_\nu^{-1}[F](x) = \frac{\pi}{2 \, {\rm sin}(\pi \, \nu)} \, 
[{\cal K}_{-\nu}^{-1}[F](x) - {\cal K}_\nu^{-1}[F](x)],
\]
for integer $\nu$, the right hand side of these relations are replaced by their limiting
values. Fortunately, Oberhettinger \cite{Oberhettinger1972book}has tabulated a rich collection of
$K$-transforms and their inverses. For more details concerning integrals
involving Bessel functions as integrands and their relationship to Laplace
transformations, see the work by Luke \cite{Luke1962book} and Erd\'elyi 
\cite{Erdelyi1954book}. 


%


\bibliographystyle{alpha}
\bibliography{Bibliog}

\end{document}